\documentclass[a4paper,12pt]{article}

\usepackage[utf8]{inputenc}
\usepackage[english]{babel}
 \usepackage{mathrsfs}
			
\usepackage{enumerate}
\usepackage{amsmath}
\usepackage{amsthm}
\usepackage{amsfonts}
\usepackage{textcomp} 
\usepackage{hyperref}
\usepackage{blindtext}
\usepackage[inline]{enumitem}

\let\oldbibliography\thebibliography
\renewcommand{\thebibliography}[1]{\oldbibliography{#1}
\setlength{\itemsep}{0pt}}

\theoremstyle{plain}
\newtheorem*{theorem*}{Theorem}
\newtheorem{theorem}{Theorem}[section]

\newtheorem{corollary}[theorem]{Corollary}
\newtheorem{proposition}[theorem]{Proposition}

\newtheorem{lemma}[theorem]{Lemma}
\newtheorem{remark}[theorem]{Remark}

\newcommand{\pr}{\mathbb{P}}

\newcommand{\E}{\mathbb{E}}
\newcommand{\R}{\mathbb{R}}
\newcommand{\Na}{\mathbb{N}}

\newcommand{\Z}{\mathbb{Z}}

\newcommand{\N}{\mathcal{N}}
\newcommand{\Lsite}{\mathbb{L}}
\newcommand{\G}{\mathscr{G}}
\newcommand{\Law}{\mathcal{L}} 

\newcommand{\B}{\mathbb{B}}
\newcommand{\cov}{\textup{Cov}}

\newcommand{\I}{\mathcal{I}} 
\newcommand{\J}{\mathcal{J}} 

\renewcommand{\L}{\mathscr{L}}

\newcommand{\Var}{\mathrm{Var}}
\newcommand{\Cov}{\mathrm{Cov}}

\title{The central limit theorem for supercritical \\ oriented percolation in two dimensions} 
\author{\textbf{Achillefs Tzioufas}\footnote{\small{\texttt{tzioufas@ime.usp.br}}}}
\begin{document}
\maketitle

\vspace{-1mm}
\begin{abstract} 
\noindent We consider the cardinality of supercritical oriented bond percolation in two dimensions.  We show that, whenever the origin is conditioned to percolate, the process appropriately normalized converges asymptotically in distribution to the standard normal law. This resolves a longstanding open problem pointed out to in several instances in the literature.  The result applies also to the continuous-time analog of the process, viz.\ the basic one-dimensional contact process.  
We also derive general random-indices central limit theorems for associated random variables as byproducts of our proof.   

\vspace{1mm}

\end{abstract}

\section{Introduction} 
We consider oriented bond percolation on the two-dimensional integer lattice. For background on this process, we refer to the review \cite{D84}. We show that the process exhibits classic central limit theorem (CLT) behavior in all of the supercritical phase; meaning that the law of the diffusively rescaled cardinality of the process started from a site conditioned to percolate converges asymptotically in distribution to the standard normal law.  The continuous-time analog of two-dimensional oriented percolation is the basic contact process in one (spatial) dimension.  Our result and approach convey analogously to this process. The contact process on integer-lattices was introduced in \cite{H74}. The corresponding strong law of large numbers (SLLN) was shown in \cite{DG83}. The CLT has been posed as an open problem originally in \cite{D84}, and later in \cite{D88}, and more recently, in \cite{D91}. The contact process falls into the subject of interacting particle systems (IPS), for background on which we refer to the classic accounts \cite{G79, L85, D88}, furthermore, we refer to the later, more comprehensive accounts  \cite{D95} and \cite{L99} as the subject received additional attention afterward in the literature. Percolation theory originates in \cite{BH57} and, for more in this regard, we refer to the classic accounts \cite{G99, BR06}.

Harris'  Growth Theorem \cite{H78} regards that the rate of growth of the highly supercritical contact process conditioned to percolate is almost surely linear in all dimensions.\footnote{see Theorem 13.5 in \cite{H78}; the appellation is due to \cite{G81}, see the Remark following Theorem 9 there, see also Theorem 1.1 in \cite{D80}.} The corresponding $L^{1}$-LLN for the supercritical process in one dimension was shown by means of subadditivity and coupling arguments in \cite{D80}.\footnote{see Theorem 1.2 in  \cite{D80}.} We note that this result  is considered a precursor to the general subadditive ergodic theorem shown in \cite{L85-2}. The range of parameter values for which the Harris' Growth Theorem holds was extended by means of improvements to Peierls' argument in continuous time, shown in \cite{G81}.\footnote{see Theorem 9 in \cite{G81}.} The SLLN for the process in all dimensions with parameter value larger than the critical value of the one-sided process in dimension one was derived, as a corollary of the general shape theorem, in \cite{DG82}. The SLLN, valid for all of the supercritical phase in dimension one, was shown by means of renormalization group techniques in \cite{DG83}.\footnote{see \cite{DG83}, Theorem 9.} Furthermore, the important property that the invariant measure possesses exponentially decaying correlations,  together with other exponential estimates, was shown there.\footnote{see \cite{DG83},  Theorems 7 and 8.}  This key property, together with the SLLN for the position of the endpoints, result shown earlier in \cite{D80},  enabled the proof of the SLLN in \cite{DG83}. Among other landmark results, the shape theorem, and hence the SLLN, valid for all of the supercritical phase and in all dimensions, was shown by means of renormalization techniques in \cite{BG90}, see also the review \cite{D91}. To date the following CLT's have been derived in the literature regarding other functionals of the contact process. The CLT regarding time-averages of finite support functions of the infinite one-dimensional supercritical contact process was shown in \cite{S86}, by means of following an approach by \cite{CG83}, using an exponential decay property by \cite{DG82}, and applying general results of \cite{NW81, N83}. Further, we mention that the CLT regarding the endpoints of the process was shown by means of mixing techniques in \cite{GP87}, and later by means of elementary arguments in \cite{K89}. In addition, for a detailed literature account regarding known CLT's in classic percolation, we refer to $\mathsection$ 11.6 in  \cite{G99}, see also the later \cite{P01, P05}. 

Furthermore, we derive certain CLT's regarding randomly-indexed partial sums of non-stationary, associated r.v.'s (random variables), as byproducts of our proof technique. To the limits of one's knowledge, randomly-indexed CLT's for families of associated r.v.'s have not been considered elsewhere in the literature. The introduction and the appreciation of the usefulness of association in percolation dates back to the Harris' Lemma \cite{H60}\footnote{see Lemma 4.1 in \cite{H60}; the appellation is attributed by \cite{BR06}.},  with most prominent extensions to non-product measures, being the FKG inequality \cite{FKG71}, the Holley inequality \cite{Ho74}\footnote{see also \cite{G12} or \cite{GHM01}.}, and the Ahlswede–Daykin inequality \cite{AD78}\footnote{see also \cite{AS16}.}. The systematic study of this concept as a general dependence structure was initiated in \cite{EPW67}. The acknowledgment of that asymptotics for the correlation structure are useful in studying approximate independence of associated r.v.'s originates in \cite{L72}, where necessary and sufficient conditions of this sort for the ergodic theorem to extend to this case were shown. The first corresponding CLT was derived in \cite{N80}, whereas the key notion of demimartingales was introduced in \cite{NW81}. Other notable CLT's, which replace the stationarity assumption with moment conditions, are those due to \cite{CG84}, and also \cite{B87}. For background and comprehensive expositions on association, classic limit theorems for associated r.v.'s in particular, and much more about recent advances on the subject, we refer to the reviews \cite{N84}, \cite{BS07}, \cite{O12}, \cite{Ra12}.

Our main result comprises the CLT for supercritical oriented percolation in two dimensions, which we state explicitly in $\mathsection$ \ref{secCLT}, Theorem \ref{CLT}. Regarding the CLT's for randomly-indexed associated r.v.'s, see  $\mathsection$ \ref{seccgprod}. 

\subsection{Definition of the process}
We let $\L = \L( \Lsite, \B)$ be the usual two-dimensional oriented percolation lattice graph,
for which the set of sites is $\Lsite = \{ (x,n) \in \Z^{2}: x+n \in 2\Z \mbox{ and } n \geq 0\}$, $2\Z = \{2k: k \in \Z\}$, and the set of bonds is $\B = \{\left[(x,n), (y, n+1)\right>: |x-y| =1\}$, and where $\left[s, u\right>$ means that an arrow (or bond) directed from site $s$ to site $u$ is present, see fig.\ 1, p.\ 1001, \cite{D84}, for this and other transpositions of $\L$ in the plane. We consider independent bond percolation on $\L$ with open (or retaining) probability parameter $p \in [0,1]$, defined as follows.  We consider the configuration space $\Omega = \{0,1\}^{\B} = \{ \omega: \B \rightarrow \{0,1\} \}$. We let $\pr (= \pr_{p})$ denote the joint distribution of $(\omega(b): b \in \B)$, an ensemble of i.i.d.\  $p$-Bernoulli r.v.'s, which is, $\mu(\omega(b)= 1) = 1- \mu(\omega(b)= 0) = p$. We note that $\pr$ yields a probability measure on $\Omega$, equipped as usual with $\mathcal{F}$, the $\sigma$-field of subsets of $\Omega$ generated by finite-dimensional cylinders.  We may further let $\G = \G( \Lsite, \B_{1})$, $\B_{1} = \{b : \omega(b) =1\}$, be the subgraph of $\L$ in which $b$ is retained if and only if $\omega(b) = 1$.  

For given $\omega \in \Omega$, bonds $b$ such that $\omega(b)= 1$ are thought of as open (or retained); whereas if $b$ is assigned value $\omega(b)= 0$, we consider $b$ as closed (or removed), which may be thought of as flow being disallowed. If  $s_{m}, s_{n} \in \Lsite$, $s_{m}= (x_{m}, m)$, $s_{n}=  (x_{n},n)$, $m \leq n$, then, given $\omega \in \Omega$, we write $ s_{m} \rightarrow s_{n}$ whenever there is a directed path from $ s_{m}$ to $s_{n}$ in $\G(\omega)$, that is, there is $s_{m+1} = (x_{m+1}, m+1), \dots, s_{n-1} = (x_{n-1}, n-1)$ such that $ \omega(\left[s_{k}, s_{k+1}\right>) = 1$ for all $m \leq k \leq n-1$. 

We let $\xi_{n}^{\eta} = \{x:  \left(y,0 \right) \rightarrow (x,n), \mbox{ for some } y \in  \eta \}$, $\eta \subseteq 2\Z$. We note also that, when convenient, we shall use the coordinate-wise notation  $\xi_{n}^{\eta}(x) = 1(x \in \xi_{n}^{\eta})$,  where we denote by $1(A)$ the indicator random variable of an event $A$. We note that $(\xi_{n}^{\eta}, \eta \subset 2\Z)$ furthermore admits a Markovian definition, and hence, we may think of the vertices' first and second coordinates as space and time, respectively. We also note that, by definition of $\Lsite$, $\xi_{n}^{\eta}\subset 2\Z$, for $n \in 2\Z_{+}$, and $\xi_{n}^{\eta}\subset 2\Z+1$, for $n \in 2\Z_{+} +1$.

We will denote simply by $(\xi_{n})$ the process started from $O = \{0\}$ and, in general, we will drop superscripts associated to the starting set in our notation when referring to $O = \{0\}$.

\subsection{The critical value and the upper-invariant measure}
We state here some basic definitions and facts, for a more detailed exposition of which, see for instance, \cite{D84, G12, L10}. 
We let $\Omega_{\infty}^{\eta}$ be the percolation event for initial configurations $\eta$, $|\eta|< \infty$, that is, we let
\begin{equation}\label{eqrhos}
\Omega_{\infty}^{\eta} := \cap_{n\geq1} \Omega_{n}^{\eta} \hspace{3mm} \mbox{ and } \hspace{3mm} \Omega_{n}^{\eta} := \{|\xi_{n}^{\eta}| \geq 1\},
\end{equation}
where $| \cdot|$ denotes cardinality; and we also note that $\Omega_{n}^{\eta} \supseteq \Omega_{n+1}^{\eta}$, $\pr$-a.s.. Further, we let $\rho = \rho(p)$ be the so-called asymptotic density, defined as follows
\begin{equation}\label{rhodef}
\rho(p) = \pr(\Omega_{\infty}) = \lim_{n \rightarrow \infty} \rho_{n}, \hspace{5mm} \rho_{n} := \pr(\Omega_{n}),  
\end{equation}
where the appellation derives in view of (\ref{definvm}), from the resultant equality of $\rho(p)$ by appropriately applying (\ref{duintro}) below.  We let in addition $p_{c}$ be the critical value,  defined as follows
\begin{equation}\label{rhocridef}
p_{c} = \inf\{ p: \rho(p)  >0\}. 
\end{equation}
Where we recall that it is elementary that $p_{c} \in (0,1)$,  and  that the well-definedness of $p_{c}$ comes from that $\rho(p)$ is non-decreasing in $p$, which is an elementary consequence of the construction by the superposition of Bernoulli r.v.'s property.
We also note that the assumption that $p>p_{c}$ we require here may be replaced by the apriori weaker assumption that $\rho(p)>0$, since the two assumptions are equivalent due to that  $\rho(p_{c})=0$,  shown for all dimensions in \cite{BG90}, see also \cite{D91}, and also \cite{BG94} for the extension of this results to general attractive spin systems. 

Let further, 
\begin{equation}\label{sigdef}
 \Sigma_{0} =\{\eta \subset 2\Z:  |\eta| < \infty \}, \hspace{4mm} \Sigma =\{\eta \subset 2\Z, |\eta| = \infty\}, 
\end{equation}  
We recall also that, if $\mu_{n}$ denotes the distribution of $\xi_{2n}^{2\Z}$, then we have that
\begin{equation}\label{definvm} 
\mu_{n} \Rightarrow\bar{\nu}, \mbox{ as  }  n \rightarrow \infty, 
\end{equation} 
where $\bar{\nu}$ is the so-called upper-invariant measure, defined on $\Sigma$ and uniquely determined by its cylinders, which,  in view of the so-called self-duality property (see, for example, (\ref{duapp}) below),  is such that
\begin{equation}\label{duintro}
\bar{\nu}(\eta: \eta \cap B \not= \emptyset) = \pr(\xi_{n}^{B} \not= \emptyset, \mbox{ for all } n \geq 1),
\end{equation}
$B \in \Sigma_{0}$, and where `$\Rightarrow$' denotes weak convergence, which we define as convergence of the finite-dimensional distributions
\[
\pr(\xi_{2n}^{2\Z} \cap B = C), \mbox{ for } C \subset B \in \Sigma_{0}, 
\]
as $n \rightarrow \infty$. To see the reason that we refer to $\rho(p)$ as the asymptotic density, note that by (\ref{duintro}) we have that $\bar{\nu}(\eta: \eta \cap O \not= \emptyset) =  \rho$.  Further, we note that (\ref{definvm}) is denoted below simply as follows,
\[
\xi_{2n}^{2\Z} \Rightarrow \bar{\xi}, \hspace{5mm} n \rightarrow \infty, 
\]
where $\bar{\xi}$ is a random field distributed according to $\bar{\nu}$, denoted as $\bar{\xi} \sim \bar{\nu}$ below.

Finally, prior to turning to our main statement in the section below,
we give some additional notation first.  
We note that the shorthands  $\Law(X_{n}) \xrightarrow{w} \mathcal{N}(0,\sigma^{2})$, as $n \rightarrow \infty$, $n \in \Na$, as well as $X_{n} \xrightarrow{w} \mathcal{N}(0,\sigma^{2})$, as $n \rightarrow \infty$, will be in force in the sequel in order to denote weak convergence to a normal distribution with mean $0$ and variance $\sigma^{2}$, which is that, if we let $F_{n}(x)$ be the cumulative distribution function associated with $X_{n}$, we have that
$F_{n}(x) \rightarrow \int_{- \infty}^{x} (2\pi)^{- 1/2} e^{- u^{2}/2} \mbox{d} u$, as $n \rightarrow \infty$.

\section{Results}

\subsection{The CLT}\label{secCLT}
To state next our main result, we let $p>p_{c}$ and we let $\bar{\xi} \sim \bar{\nu}$. We let also 
\begin{equation}\label{eqsig}
\sigma^{2}  =  \sum_{x} \cov(x \in\bar{\xi},O \in \bar{\xi}) < \infty.    
\end{equation}
Furthermore, we let $\bar{\pr}$ be the probability measure induced by the original $\pr$ by conditioning on $\Omega_{\infty}$, which is, $\bar{\pr}(\cdot) = \pr(\cdot|  \Omega_{\infty})$.  We denote by $\Law(X| \Omega_{\infty})$ the law of a r.v.\ $X$ under $\bar{\pr}$.  
In addition, we let $r_{n} = \sup\xi_{n}$ and $l_{n} = \inf\xi_{n}$, and we further let $d_{n} = \frac{1}{2} (r_{n} - l_{n}) +1$. We recall also that $\rho_{n} = \pr(\Omega_{n})$. 

\begin{theorem}\label{CLT} Let $p > p_{c}$. We have that, as $n \rightarrow \infty$, 
\begin{equation*}\label{CLT}	
\Law\left(\frac{ | \xi_{n}| - d_{n} \rho_{n} } {\sigma \sqrt{d_{n}}} \mbox{ } \vline \mbox{ } \Omega_{\infty} \right) \overset{w}{\longrightarrow} \N(0, 1). 
\end{equation*}
\end{theorem} 

From a technical perspective, the main novelty in our proof approach is the consideration of the (non-stopping) time regarding the last-intersection of the two infinite endpoints processes.  To briefly elaborate on this coupling observation here (see the Lemmas subsequent to its definition in (\ref{taf}) for the exact statements), we let 
\[
r^{-}_{n} = \sup \xi_{n}^{2\Z_{-}}  \mbox{ and } l_{n}^{+}= \inf \xi_{n}^{2\Z_{+}}, 
\] 
where $2\Z_{-} = \{ \dots,-2, 0\}$, and $2\Z_{+} = \{ 0, 2, \dots \}$. We note that, as will be seen from the proof, the distribution of $|\xi_{n}| = \sum_{x=l_{n}}^{r_{n}} \xi_{n}(x)$ conditioned on $\Omega_{\infty}$, is equal to that of $\sum_{x=l_{n}^{+}}^{r^{-}_{n}} \xi_{n}^{2 \Z}(x)$, 
for all $n$ after this random time occurs, and therefore, asymptotics for the latter process permit to infer the same asymptotics for the former one. In this manner, we circumvent the effects of altering the distribution of $\xi_{n}$ when conditioning on $\Omega_{\infty}$, and thus, we are able to deduce Theorem \ref{CLT} by working on the whole probability space, and dealing with partial sums of the infinite processes, involved in $\sum_{x=l_{n}^{+}}^{r^{-}_{n}} \xi_{n}^{2 \Z}(x)$, instead. We further note that, in order to deal with the fact that this is a non-stopping time in our proof, we reside on independence inherited from the independence of the underlying Bernoulli r.v.'s in disjoint parts of $\L$, an observation applied in a different context by \cite{K89}. We also note that this coupling is intrinsic to two-dimensions, since it relies on path intersection properties, and that hence, we expect that new methods will be required for the extension of this result to higher dimensions. 
On the other hand, we believe and pursue in forthcoming work \cite{Tfor} that the techniques we develop can be used in order to give a proof of the law of the iterated logarithm corresponding to Theorem \ref{CLT}. Furthermore, we note that the method of proof of Theorem \ref{CLT} relies on Proposition \ref{AnscC}, stated in $\mathsection$ \ref{secans} below, and also incorporates an earlier observation due to \cite{D84}.  We note that a key ingredient for Proposition \ref{AnscC} to apply in the context of supercritical oriented percolation is the Harris' correlation inequality \cite{H77}; see also Theorem B.17 in \cite{L99} and the references therein.  In addition, we note that our proof approach, and the techniques involved, differentiate from those devised in known CLT's for percolation processes, due to the fact that we consider partial sums that are indexed randomly, depending on the state of the process itself.

\subsection{Random-indices CLT's}\label{seccgprod}
We recall here the following definition.

{\bf Association. } A collection of r.v.'s $(X_{i}: i \in I)$, $|I| = \infty$, is associated if for all finite sub-collections $X_{1}, \dots, X_{m}$ and all coordinate-wise non-decreasing $f_{1}, f_{2}: \R^{m} \rightarrow \R$ we have that $\Cov(\tilde{f}_{1}, \tilde{f}_{2}) \geq 0$, $\tilde{f}_{j} := f_{j}(X_{1}, \dots, X_{m})$, $j=1,2$, whenever this covariance exists.

By known results our proof approach 
provides with certain random-indices central limit theorems for associated triangular arrays of r.v.'s, which we effectively obtain as direct byproducts. 
One important aspect of those statements is that nothing is assumed regarding independence among the summands and the index family of r.v.'s. Corollary \ref{cgprod} stated next in particular is a random-index extended version of the CLT in Theorem 1 due to \cite{CG84} with the additional proviso (\ref{ansacor}) below. Furthermore, we note that 
we may in addition obtain in a manner which is directly analogous and is thus omitted the corresponding random-index CLT's extensions to Theorem 3 in \cite{B87}, or Theorem 3 in \cite{NW81}, with the said additional proviso. 

\begin{corollary}\label{cgprod}
Let $\{X_{n}(j): 0 \leq j \leq n\}$ be such that 
$\E(X_{n}(j))= 0$, $\forall$ $n, j$, and that, for each $n$,
\begin{equation}\label{eqcltassoc}
\{X_{n}(j)\} \mbox{ are associated}. 
\end{equation}
Suppose also that 
\begin{equation}\label{cg1}
\inf_{j, n} \Var(X_{n}(j))>0 \hspace{3mm}\mbox{ and }\hspace{3mm} \sup_{j, n}\E(|X_{n}(j)|^{3}) < \infty. 
\end{equation}
Furthermore, suppose that $u(r) = \sup_{j, n}\sum_{|k - j| \geq r} \Cov(X_{n}(j), X_{n}(k))$, $r \geq0$, is such that
\begin{equation}\label{cg2}
 u(r) < \infty, \mbox{ for all } r,  \mbox{ and that } u(r) \rightarrow 0, \mbox{ as } r \rightarrow \infty.
\end{equation}
Let $S_{n}(i) = \sum_{j=0}^{i} X_{n}(j)$, and assume in addition that  
\begin{equation}\label{ansacor}
\sup_{j, n}\Cov\left(X_{n}(j), S_{n}(j-1) \right) <\infty. 
\end{equation}
Let $(N_{n}, n \in \Na)$ be integer-valued and positive r.v.'s, such that  
\begin{equation}\label{Nclt}
\frac{N_{n}}{n} \xrightarrow{w} \theta,  \mbox{ as }  n \rightarrow \infty, 
\end{equation}
for some $0<\theta \leq 1$. 

We then have that 
\[
\frac{S_{n}(N_{n})}{\sqrt{N_{n}}} \xrightarrow{w} \mathcal{N}(0,\sigma^{2}),  \mbox{ as } n \rightarrow \infty
\]
and also that 
\[
\frac{S_{n}(N_{n})}{\sqrt{ \theta n}} \xrightarrow{w} \mathcal{N}(0,\sigma^{2}),  \mbox{ as } n \rightarrow \infty,
\]
where $\sigma^{2}:= \lim_{n \rightarrow \infty} \textup{\mbox{Var}} (S_{[ \theta n]} / \sqrt{[\theta n]})$, $0<\sigma^{2}<\infty$. 
\end{corollary} 

By Theorem 1 in \cite{CG84}, the method of proof of Corollary \ref{cgprod} further relies on an application of Lemma \ref{Ansc}, which we derive on the way to the proof of Proposition \ref{AnscC} below.

\subsection{Anscombe's condition}\label{secans}
Proposition \ref{AnscC} next regards a condition about deviations of random partial sums from deterministic ones the interval. This condition in the i.i.d.\ case was shown in \cite{A52}.\footnote{Hence the appellation of the condition attributed to by \cite{GU09}.}  The validity of this condition has not been anticipated to extend in the generality of Proposition \ref{AnscC}, see Remark \ref{heur}. To state it, we write $X_{n} \xrightarrow{p} X$, as $n \rightarrow \infty$, to denote convergence in probability. Further, we let $\{X_{t}(j):  (j, t) \in \L \}$ and let $S_{t}(u, v) = \sum_{j =u}^{v} X_{t}(j)$. We introduce the following assumptions, which we will invoke there. 
\begin{equation}\label{mnzero}
\E(X_{t}(j))= 0, \mbox{ for all } (j, t) \in \Lsite,
\end{equation}  
\begin{equation}\label{Cassass}
\{X_{t}(j)\} \mbox{ is associated for each } t, 
\end{equation}
\begin{equation}\label{Cans0}
\sup_{j, t}\E(X_{t}(j)^{2})< \infty;
\end{equation}
furthermore, let $S_{t}^{+}(v) = \sum_{j = 0}^{v} X_{t}(j)$,  $S_{t}^{-}(u) = \sum_{j = 0}^{u} X'_{t}(j)$, $X'_{t}(j) = X_{t}(-j-1)$, $u, v \geq 0$,  $j\geq 0$ and assume that
\begin{equation}\label{Cansa}
C^{+} = \sup_{j, t \geq 0}\Cov\left(X_{t}(j), S_{t}^{+}(j-1) \right) <\infty, C^{-} := \sup_{j, t \geq 0}\Cov\left(X_{t}'(j), S_{t}^{-}(j-1) \right) <\infty. 
\end{equation}
We further let $(M_{t}: t\geq 0)$ and $(m_{t}: t\geq 0)$ be such that  $(M_{t}, t) \in \Lsite$ and that $(m_{t}, t) \in \Lsite$; we assume that, for some $0<\theta <\infty$,
\begin{equation}\label{Cansb}
\frac{M_{t}}{t} \xrightarrow{w} \theta \hspace{3mm} \mbox{ and } \hspace{3mm}  \frac{m_{t}}{t} \xrightarrow{w} - \theta, \mbox{ as } t \rightarrow \infty. 
\end{equation}

\begin{proposition}\label{AnscC}
We let $\{X_{t}(j):  (j, t) \in \Lsite \}$  and let $S_{t}(u, v) = \sum_{j =u}^{v} X_{t}(j)$. Let us assume that conditions (\ref{mnzero}), (\ref{Cassass}),  (\ref{Cans0}), and (\ref{Cansa}) are fulfilled. 
We let $(M_{t}: t\geq 0)$ and $(m_{t}: t\geq 0)$ be such that  $(M_{t}, t) \in \Lsite$ and that $(m_{t}, t) \in \Lsite$ and, further, assume that (\ref{Cansb}) is fulfilled. We then have that 
\begin{equation}\label{SSzero}
\frac{S_{t}(m_t, M_{t}) - S_{t}(-\theta t, \theta t)}{ \sqrt{ \theta t} } \xrightarrow{p} 0,  \mbox{ as } t \rightarrow \infty.
\end{equation} 
\end{proposition}

\noindent  Where we note that throughout here, and in the above statement in particular, we will write that  
$\sum_{x = -c}^{C}$ for $\sum_{x = -[c] -1 }^{ [C] }$, where $[\cdot]$ denotes the largest integer smaller than the argument, and that we also use the notational convention $\sum_{0}^{-1} := 0$.

The method of proof of Proposition \ref{AnscC} extends the direct proof approach due to \cite{R60} for showing the Anscombe condition in the case of i.i.d. summands.  We note that our proof invokes the so-called Hajek-R\' enyi inequality for associated r.v.'s, due to \cite{C00}. The proof of Theorem \ref{CLT} given here relies on Proposition \ref{AnscC} and thus, follows an elementary approach, see also Remark \ref{another} for a different approach. The random-index CLT's in $\mathsection$ \ref{seccgprod}, as we noted above, are in addition consequences of Proposition \ref{AnscC}, which we find of independent interest.

{ \bf Outline of Proofs. } The remainder of this paper is organized as follows. The proof of Theorem \ref{CLT}, by means of applying Proposition \ref{AnscC}, is given in $\mathsection$ \ref{secthm1}. Preliminaries we will invoke in this proof are stated first in $\mathsection \mathsection$  \ref{bas} separately, whereas another proof of Proposition \ref{bypro} stated below in there is provided with for completeness in the Appendix  $\mathsection$ \ref{secApp}. In  $\mathsection$ \ref{randev}, the proof of Proposition \ref{AnscC} is provided with, see $\mathsection \mathsection$ \ref{subspro}. That of Corollary \ref{cgprod} is also given there, in $\mathsection \mathsection$ \ref{prfcor}.

\section{Theorem \ref{CLT}}\label{secthm1}

\subsection{Preliminaries}\label{bas}

We briefly state certain facts on oriented percolation that we use later on.

{ \bf Some notation. }The following definitions, which will also be useful for simplifying notations below, are introduced.  Recall that we let $r^{-}_{n} = \sup \xi_{n}^{2\Z_{-}}$ and $l_{n}^{+}= \inf \xi_{n}^{2\Z_{+}}$, where $2\Z_{-} = \{ \dots,-2, 0\}$, and $2\Z_{+} = \{ 0, 2, \dots \}$.
We let
\begin{equation}\label{iotacap}
\I_{n} = \{x: l_{n} \leq x \leq r_{n}, (x,n) \in \Lsite\} 
\end{equation}
if $l_{n} \leq r_{n}$, and $\I_{n} = O$, otherwise.  Similarly, we let 
\begin{equation}\label{jotacap}
\J_{n} = \{x: l_{n}^{+} \leq x \leq  r_{n}^{-}, (x,n) \in \Lsite\}
\end{equation}
if $l_{n}^{+} \leq r_{n}^{-}$, and $\J_{n} = O$, otherwise. To see our motivation for considering $\J_{n}$, and $\I_{n}$ analogously, note that
\begin{equation*}
|\J_{n}| = \frac{r_{n}^{-} - l^{+}_{n}}{2} + 1, \hspace{1mm} \mbox{ on }\hspace{1mm} \{l^{+}_{n} \leq r_{n}^{-}\}, 
\end{equation*}
and $|\J_{n}| = 1$, otherwise.   We let in addition the family of centered r.v.'s, which will play a central r\^ole in our analysis below, $(\hat{\xi}_{n}^{2\Z}(x): x \in 2\Z)$,  as follows.  We let
\begin{equation}\label{centered}
\hat{\xi}_{n}^{2\Z}(x) = \xi_{n}^{2\Z}(x) - \rho_{n}, \mbox{ for all } n \geq 1,  
\end{equation}
where $\rho_{n} = \pr(\Omega_{n})$ is defined in  (\ref{rhodef}).  We note that $(\hat{\xi}_{n}^{2\Z}(x): x \in 2\Z)$ are zero-mean, since by (\ref{duapp}) below, $\E(\xi_{n}^{2\Z}(x)) =  \rho_{n}$.

{\bf The basic coupling. }
We state an important observation due to \cite{D80}, which comprises the following consequence of path intersection properties.  We have that 
\begin{equation}\label{basiccoup}
\xi_{n} =  \xi_{n}^{2 \Z} \cap [l_{n}, r_{n}]  = \xi_{n}^{2 \Z} \cap [l_{n}^{+}, r_{n}^{-}] 
\hspace{5mm} \mbox{ on } \Omega_{n}, 
\end{equation}
$\pr \mbox{-a.s.}$ and, in particular, 
\begin{equation}\label{couprl}
r_{n} =  r_{n}^{-} \hspace{3mm}  \mbox{ and } \hspace{3mm}  l_{n} = l_{n}^{+}, \hspace{5mm} \mbox{ on } \Omega_{n}, 
\end{equation}
$\pr \mbox{-a.s.}$ and further, (\ref{basiccoup}) gives that
\begin{equation}\label{cardi}
|\xi_{n}| = \sum_{x \in \J_{n}} \xi_{n}^{2\Z}(x),  \mbox{ }  \hspace{5mm}  \mbox{ on }  \Omega_{n}
\end{equation}
$\pr \mbox{-a.s.}$. Furthermore, since $\Omega_{n} = \{r_{k} \geq l_{k}, \forall k \leq n\}$, we also have that 
\begin{equation}\label{omeend}
\Omega_{n} = \{r_{k}^{-} \geq l_{k}^{+}, \forall k \leq n \}.
\end{equation}
Further, we note that (\ref{omeend}) is in fact a special case of the following statement, regarding general initial configurations.  Let $\eta^{-}$ and $\eta^{+}$ be such that $\eta^{-}(O) =\eta^{+}(O) =1$, and also $\eta^{-}(x)=0$, for $ x\geq 2$, whereas, $\eta^{+}(x) = 0$, $x \leq -2$, and otherwise arbitrary. Letting $r_{n}^{\eta^{-}} = \sup\{x:  \xi_{n}^{\eta^{-}}(x) =1 \}$ and $l_{n}^{\eta^{+}}= \inf\{x: \xi_{n}^{\eta^{+}}(x) = 1\}$, we have that, for all $n \geq 1$, on $\Omega_{n}$,  $r_{n}  = r^{\eta^{-}}_{n}$ 
and $l_{n} = l_{n}^{\eta^{+}}$ and, further that
\begin{equation}\label{rlc}
\Omega_{n} = \{r_{m}^{\eta^{-}} \geq  l_{m}^{\eta^{+}}, \mbox{ for all } m \leq n\}. 
\end{equation}
For proofs of these statements, see for instance, $\mathsection$3, \cite{D84}, see also \cite{D80, G81}.

{ \bf The asymptotic velocity. }
For all $p > p_{c}$, there is $\alpha = \alpha(p)>0$, such that 
\begin{equation}\label{avel}
\lim_{n \rightarrow \infty} \frac{r^{-}_{n}}{n} = \alpha \hspace{2mm} \mbox{ and }  \hspace{3mm} \lim_{n \rightarrow \infty}\frac{l^{+}_{n}}{n} = -\alpha, 
\end{equation}
$\pr$-a.s.. Further, we have that (\ref{avel}) yields from (\ref{couprl}), that 
\begin{equation}\label{spd}
\lim_{n \rightarrow \infty} \frac{r_{n}}{n} = \lim_{n \rightarrow \infty} \frac{-l_{n}}{n} = \alpha, \mbox{ } \bar{\pr}\mbox{ a.s.}, 
\end{equation}
where we refer to $\alpha := \alpha(p)$ as the asymptotic velocity. For a proof of (\ref{avel}) we refer to Theorem 1.4 in \cite{D80}, and also (7) in $\mathsection$ 3 in \cite{D84}. 

{ \bf The SLLN.} Let $p > p_{c}$.  Let $\rho$ and $\alpha$ be the asymptotic density and velocity, as defined in (\ref{avel}) and in (\ref{rhodef}), respectively. We have that 
\begin{equation}\label{sldg}
\lim_{n \rightarrow \infty} \frac{|\xi_{n}| }{n} = \alpha \rho, \mbox{ } \bar{\pr} \mbox{-a.s.}, 
\end{equation}
For a proof of (\ref{sldg}) we refer to Theorem 9 in \cite{DG83}, see also (2) in $\mathsection$ 13 in \cite{D84}.  We mention here that from (\ref{basiccoup}) and (\ref{spd}), since $|\I_{n}| = \frac{r_{n}-l_{n}}{2}+1$, we have that, as $n \rightarrow \infty$, $\displaystyle{\frac{|\I_{n}|}{n} \rightarrow \alpha \mbox{, } \bar{\pr}\mbox{ a.s.}}$. Thus, we have that (\ref{sldg}) yields that, as $n \rightarrow \infty$, $\displaystyle{\frac{\sum_{x \in \I_{n}}  \xi_{n}(x)}{|\I_{n}|} \rightarrow \rho, \hspace{2mm}  \bar{\pr}\mbox{ a.s.}}$.

\textbf{Large deviations.} Let $p > p_{c}$. let $a< \alpha(p)$. Then, the following limit exists and is strictly negative, 
\begin{equation}\label{ldr}
\lim_{n \rightarrow \infty} \frac{1}{n} \log \pr(r_{n}^{-} < a). 
\end{equation}  
We require in addition below, the following known elementary consequence of (\ref{ldr}). There are $C, \gamma \in (0, \infty)$, such that 
\begin{equation}\label{expoldev}
\pr(\exists m\geq n: r^{-}_{m} < 0)\leq Ce^{- \gamma n},
\end{equation}
$n \geq 1$, see, for instance \cite{D84}, p.\ 1031, $\mathsection$ 12, first display in the proof of (1). 

{\bf Monotonicity and Self-duality.} An immediate consequence of the construction is that 
\begin{equation}\label{eqmon}
A \subseteq B \hspace{2mm} \Longrightarrow \hspace{2mm} \xi_{n}^{A} \subset \xi_{n}^{A}, 
\end{equation}
$\pr$-a.s. Further, we have that, for all $n$ even,
\begin{equation}\label{eqduality}
\pr(\xi_{n}^{A} \cap B \not= \emptyset) = \pr(\xi_{n}^{B} \cap A \not= \emptyset),
\end{equation}
$A, B \subset 2\Z$, and analogously, for $n$ odd. The proof of (\ref{eqduality}), see $\mathsection$ 8, (2), p.\ 1021 in \cite{D84} relies on the observation that, after reversing the direction of all arrows in any realization, the law of the process started from $(B,2n)$, defined analogously by these new paths, and going backwards in time is the same as that of $(B,0)$; and, moreover, that a path connecting $(A,0)$ to $(B, 2n)$ exists in the original sample point if and only if there is a backwards in time path connecting $(B, 2n)$ to $(A,0)$ in the corresponding sample point.  By an application of (\ref{eqduality}) and by the definition of  the upper invariant measure  $\bar{\nu}$, see (\ref{definvm}), we note that, 
\begin{eqnarray}\label{duapp}
\pr(\xi_{n}^{A}  \not=  \emptyset, \mbox{ for all } n \geq 1 ) &= & \lim_{n \rightarrow \infty} \pr(\xi_{2n}^{2\Z} \cap A  \not= \emptyset) \nonumber \\ & = & \bar{\nu}( \eta: \eta \cap A \not= \emptyset), 
\end{eqnarray}
$A \subset \Sigma_{0}$. Furthermore, by (\ref{duapp}) and recalling the definition of $\rho$ from (\ref{rhodef}), gives that 
\begin{equation}\label{duapp2}
\rho  = \bar{\nu}( \eta: \eta \cap \{O\} \not= \emptyset) = \lim_{n \rightarrow \infty} \E(\xi_{2n}^{2 \Z}(O)). 
\end{equation}

{\bf CLT for the upper invariant measure: decay of correlations.} 
Whenever $p > p_{c}$, the upper  invariant measure $\bar{\nu}$ possesses positive, and exponentially decaying, correlations. That is, if $\bar{\xi} \sim \bar{\nu}$, we have that there are $C, \gamma \in (0, \infty)$, such that 
\begin{equation}\label{expon}
0 \leq \Cov(\bar{\xi}(0), \bar{\xi}(x)) \leq C e^{- \gamma x}, 
\end{equation}
$x \in 2\Z$. As pointed out to in \cite{DG83}, p.\ 2, see also the final Remark in $\mathsection$ 6 in \cite{G81}, property (\ref{expon}) implies the following Lemma by general results for random fields, see for instance, Theorem 12 in \cite{N84}, or the list of references before statement Proposition 4.18 in Chpt.\ I, \cite{L85}.  
\begin{lemma}\label{lem0} $\displaystyle{\Law\left(\frac{\sum_{x = -\alpha n}^{\alpha n} (\bar{\xi}(x) - \rho)}{\sqrt{\alpha n}}\right)
 \xrightarrow{w} \mathcal{N}(0, \sigma^{2})}$, as $n \rightarrow \infty$, $\sigma^{2} < \infty$,
\end{lemma} 
\noindent where $\sigma^{2} < \infty$ because $\bar{\xi}$ is strictly stationary (translation invariant), and we thus have that $\sigma^{2} = \Var(\bar{\xi}(0)) + 2 \sum_{x \geq 2} \Cov(\bar{\xi}(0), \bar{\xi}(x))$, and $\Var(\bar{\xi}(0)) = \rho - \rho^{2}$, so that $\sigma^{2} < \infty$, by (\ref{expon}). 

Furthermore in \cite{DG83} the following stronger than (\ref{expon}) property is shown.  To state it, consider $(\hat{\xi}_{n}^{2\Z}(x): x \in 2\Z)$, as defined in (\ref{centered}). We have that, for all $p >p_{c}$, there are $C, \gamma \in (0, \infty)$, such that, for any $n$, and $(x_{i} \in 2\Z: i =1, \dots, k)$, $k < \infty$, $|x_{i}- x_{j}| > 2m$, we have that  
\begin{equation}\label{expdeccor}
\left|\E\left( \prod_{i = 1}^{k} \hat{\xi}_{n}^{2\Z}(x_{i})\right) \right| \leq Ce^{-\gamma m}, 
\end{equation}
and we refer to Theorem 8 in \cite{DG83}, see also (1), p.\ 1033, \cite{D84}, for a proof of (\ref{expdeccor}). We also finally mention two other routes to derive Lemma \ref{lem0}. One of them is provided with in the discussion prior to Theorem 3.23 in Chpt.\ VI, \cite{L85}. This approach relies on deriving, by means of (\ref{ldr}),  that the convergence in (\ref{definvm}) occurs exponentially fast,  which then implies as shown there by general results that $\bar{\nu}$ has exponentially decaying correlations, from which the conclusion follows as noted above. The other route is provided with in $\mathsection$ 6 of \cite{G81}, where Lemma \ref{lem0} is derived under the condition that, there exists $C, \gamma \in (0, \infty)$,  such that 
$\pr(\bar{\Omega}_{\infty}^{\{0, \dots, 2n\}}) \leq C e ^{-\gamma n},$ for all $n \geq 1$, which is shown there to be valid for sufficiently large values of $p$, and later shown for all $p> p_{c}$ in \cite{DG83}, where we note that $\bar{E}$ denotes the complement of event $E$.

{\bf CLT for the infinite process in a cone.} We state here an observation, pointed out in \cite{D84} see $\mathsection$ 13, (4); see also p.\ 286 in \cite{D88}. Property (\ref{expdeccor}), together with the corresponding extension of Theorem 3 in \cite{NW81} to triangular arrays, yields that $\xi_{n}^{2 \Z}$ obeys classic CLT behavior. To state this explicitly, recall the definition of $\hat{\xi}_{n}^{2\Z}(x) = \xi_{n}^{2\Z}(x) - \rho_{n}$ in (\ref{centered}). We let $p >p_{c}$ and let $\alpha>0$ be the associated asymptotic velocity, and $\bar{\xi} \sim \bar{\nu}$, and $\sigma^{2}:= \sum_{x \in 2\Z} \Cov(\bar{\xi}(O), \bar{\xi}(x))$, as in (\ref{eqsig}). 
\begin{proposition}\label{bypro}  
$\displaystyle{\mathcal{L}\left(\frac{\sum_{x = - \alpha n}^{\alpha n} \hat{\xi}_{n}^{2\Z}(x)}{\sqrt{\alpha n}}\right) \xrightarrow{w} \mathcal{N}(0, \sigma^{2}),\mbox{ as } n \rightarrow \infty,}$ and $\sigma^{2} < \infty$.
\end{proposition}

\begin{remark}\label{heur}
\textup{The following heuristic, which is suggested  from properties of the basic coupling above, is stated next in (5) there: 
\[
|\xi_{n}^{2\Z} \cap  [l_{n}^{+}, r_{n}^{-}]\| - |\xi_{n}^{2\Z} \cap  [-\alpha n, \alpha n]\| \approx  \rho_{n}\frac{r_{n}^{-}-\alpha n}{2} + \rho_{n}\frac{+\alpha n -l^{+}_{n}}{2} 
\]
Note that as expected there, and proved later in \cite{GP87}, and in \cite{K89}, when diffusively normalized each term of the RHS converges asymptotically to a normal distribution, which then suggests that the RHS fluctuations diffusively rescaled would not converge in distribution to zero; see also the form of the variance conjectured in the latter reference for Theorem \ref{CLT}. We note that Proposition \ref{AnscC} will allow us to show that when normalized, the LHS converges in law to zero.}
\end{remark}

{\bf An exponential estimate. }  We require below the following known estimate. Recall that $\rho_{n} =\pr(\Omega_{n})$ $\rho = \pr(\Omega_{\infty})$.  Let $p > p_{c}$. There are $C, \gamma \in (0, \infty)$, such that 
\begin{equation}\label{lemexp} |\rho_{n} - \rho|  \leq  C e^{- \gamma n},   
\end{equation} 
$n \geq 1$. To see that (\ref{lemexp}) follows from known facts note that, if $p > p_{c}$, then there are $C, \gamma \in (0, \infty)$, such that 
\begin{equation*}\label{eqlifet}
\pr(\Omega_{n} \cap \bar{\Omega}_{\infty}) \leq  C e^{- \gamma n}, 
\end{equation*}  
$n \geq 1$, where for a proof of the above display, see \cite{DG83}, see also [(1), p.\ 1031, \cite{D84}]. By the law of total probability, and because $\Omega_{n} \supseteq  \Omega_{\infty}$, we have that 
\[
\pr(\Omega_{n}) - \pr(\Omega_{\infty}) = \pr(\Omega_{n} \cap \bar{\Omega}_{\infty}). 
\]
Combining the two displays above, and noting that by definition, if $m \leq n$, then $\Omega_{n} \subseteq \Omega_{m}$, and therefore $\pr(\Omega_{n}) - \pr(\Omega_{\infty}) \geq 0$, we arrive at (\ref{lemexp}).

{ \bf Elementary facts.} We give next certain elementary probability statements. To this end,  we let $(X_{n}: n\geq 0)$ and $(Y_{n}: n \geq 0)$ be collections of r.v.'s.  For a proof of Lemma \ref{contog} stated next, see for instance, 5.11.4 in \cite{GU12}. Lemma \ref{strimpwk} regards the the basic fact that almost sure convergence is stronger than convergence in distribution, see for instance, Theorem 5.3.1 in \cite{GU12}. Finally, Lemma \ref{preltaf} follows by noting that, as $n \rightarrow \infty$, $X_{n}(\omega) - X_{k + n}(\omega) \rightarrow 0$, $\forall \omega \in \{\tau = k\}$ and any $k \geq 0$, and then considering the partition $\cup_{k \geq 0} \{\tau = k\}$, to conclude that $X_{n} - X_{\tau + n} \rightarrow 0$ a.s..

\begin{lemma}\label{contog} We have that
\begin{equation}\label{cteq1}
\Law(X_{n}) \xrightarrow{w} X \mbox{ and } \Law(X_{n} - Y_{n}) \xrightarrow{w}  0   \Longrightarrow \Law(Y_{n}) \xrightarrow{w} X, 
\end{equation}
as $n \rightarrow \infty$. Furthermore, 
\begin{equation}\label{sluts}
\Law(X_{n}) \xrightarrow{w} \gamma \mbox{ and } \Law(Y_{n}) \xrightarrow{w} Y   \Longrightarrow  \Law\left(\frac{Y_{n}}{X_{n}}\right) \xrightarrow{w} \frac{Y}{\gamma}, 
\end{equation}
$\gamma \in \R \backslash \{0\}$, as $n \rightarrow \infty$.
\end{lemma}

\begin{lemma}\label{strimpwk}
We have that, as $n \rightarrow \infty$,
\begin{equation}
X_{n} \rightarrow X \mbox{ a.s.} \Longrightarrow \Law(X_{n}) \xrightarrow{w} X.   
\end{equation}
\end{lemma}

\begin{lemma}\label{preltaf}
Let $\tau< \infty$ a.s., but otherwise arbitrary. We have that $X_{\tau + n} - X_{n} \rightarrow 0$ a.s.. 
\end{lemma}

\subsection{Proof of Theorem \ref{CLT}}\label{prfthm1}

The contents of this section comprise primarily the proof of Theorem \ref{CLT} and may be outlined as follows. This proof is given here first, by means of relying on Propositions \ref{eqrlmp} and \ref{lemfindis} we state in it and prove immediately afterward in the same order as stated. Proofs of various auxiliary statements, denominated Lemmas, required along the way in our proofs are further postponed, in order not to interrupt its course, to the end of this section. In Remark \ref{another} We discuss here briefly in regard to modifications required to obtain the IP's associated to our Theorem \ref{CLT}
\

Certain remarks regarding other route approaches are provided with at the end of these proofs. 

\begin{proof}[Proof of Theorem  \ref{CLT}.] 
Note that, due to that $\I_{n} = \{x: l_{n} \leq x \leq r_{n}, (x,n) \in \Lsite\}$, the result comprises the statement
\begin{equation*}
\mathcal{L} \left.\left(\frac{\sum_{x \in \I_{n}}  \xi_{n}(x) -  |\I_{n}| \rho_{n}}{ \sqrt{|\I_{n}|}} \right\vert \Omega_{\infty}\right)  \xrightarrow{w}  \mathcal{N}(0, \sigma^{2}),\mbox{ as } n \rightarrow \infty. 
\end{equation*}  
Taking into account also that $\sum_{x \in \I_{n}}( \xi_{n}(x) - \rho_{n}) = |\xi_{n}| -  | \I_{n}|\rho_{n}, \mbox{ on } \Omega_{\infty},$ we have that, if we let $A_{n}  = \frac{\sum_{x \in \I_{n}}(\xi_{n}^{O}(x) - \rho_{n})}{\sqrt{|\I_{n}|}}$, then to complete this proof, it suffices to show that
\begin{equation}\label{AnOme}
\mathcal{L}(A_{n}| \Omega_{\infty}) \xrightarrow{w} \N(0, \sigma^{2}), \mbox{ as } n \rightarrow \infty. 
\end{equation}

We state next a key Proposition and subsequently state an auxiliary Lemma we require. We recall that we provide with the proofs of Propositions and Lemmas stated here afterward.   
Recall first that $r^{-}_{n} = \sup \xi_{n}^{2\Z_{-}}$ and $l_{n}^{+}= \inf \xi_{n}^{2 \Z^{+}}$, $2\Z_{-}= \{ \dots,-2, 0\}$, $2 \Z^{+} = \{ 0, 2, \dots \}$. 
Recall further that we let $(\hat{\xi}_{n}^{2\Z}(x): (x,n) \in \Lsite)$, $\hat{\xi}_{n}^{2\Z}(x) = \xi_{n}^{2\Z}(x) - \rho_{n}$,   $\rho_{n} = \pr(\Omega_{n})$, as defined in (\ref{centered}). Recall also that 
$\J_{n} = \{x: l_{n}^{+} \leq x \leq  r_{n}^{-}, (x,n) \in \Lsite\}$, whenever $l_{n}^{+} \leq r_{n}^{-}$, and $\J_{n} = O$, otherwise, as in (\ref{jotacap}).

\begin{proposition}\label{eqrlmp}
Let $\displaystyle{\bar{A}_{n} = \frac{\sum_{x \in \J_{n}} \hat{\xi}_{n}^{2\Z}(x)}{\sqrt{|\J_{n}|}}}$. We have that 
 \begin{equation}\label{eq12}
\mathcal{L}(\bar{A}_{n}) \xrightarrow{w} \mathcal{N}(0, \sigma^{2}),\mbox{ as } n \rightarrow \infty. 
\end{equation}
\end{proposition}
\noindent We state the other key Proposition below, after the following auxiliary Lemma required first.  

\begin{lemma}\label{lemtaf}
Let 
\begin{equation}\label{taf}
\tau := \inf\{n \geq0: {r^{-}_{n}= l^{+}_{n}} \mbox{ and  } {r^{-}_{m} \geq l^{+}_{m}} \mbox{ } \forall \mbox{ }m> n \}. 
\end{equation}
We have that $\tau< \infty$, a.s.
\end{lemma}
\noindent We may now give the said Proposition. 
\begin{proposition}\label{lemfindis} 
$\mathcal{L}(\bar{A}_{n + \tau}) = \mathcal{L}(A_{n}| \Omega_{\infty})$, for all $n \geq 0$. 
\end{proposition}

\noindent Note that, Lemma \ref{lemtaf} by an application of Lemma \ref{preltaf} gives that, as $n \rightarrow \infty$,
\begin{equation}\label{eqas0}
\bar{A}_{n + \tau} - \bar{A}_{n} \rightarrow 0, \mbox{ a.s.}.
\end{equation}
However (\ref{eqas0}) by Lemma \ref{strimpwk} gives that $\Law(\bar{A}_{n + \tau} - \bar{A}_{n}) \rightarrow 0$, as $n \rightarrow \infty$, and hence, Proposition \ref{eqrlmp} by an application of (\ref{cteq1})an application of (\ref{cteq1}) yields that 
\begin{equation}\label{normshoft}
\mathcal{L}(\bar{A}_{n + \tau}) \xrightarrow{w} \N(0, \sigma^{2}), \mbox{ as } n \rightarrow \infty.
\end{equation} 
\noindent From Proposition \ref{lemfindis} and (\ref{normshoft}), we have that (\ref{AnOme}) follows,  and, therefore, the proof is complete.
\end{proof}

\begin{proof}[Proof of Proposition \ref{eqrlmp}.] 
From Proposition \ref{bypro} we have that 
\[
\mathcal{L}\left(\frac{\sum_{x = - \alpha n}^{ \alpha n} \hat{\xi}_{n}^{2\Z}(x)}{\sqrt{\alpha n}}\right)  \xrightarrow{w} \mathcal{N}(0, \sigma^{2}),\mbox{ as } n \rightarrow \infty
\] 
where we recall that $\alpha: = \alpha(p)>0$, $p > p_{c}$, is the asymptotic velocity, as defined in (\ref{avel}); however, note that (\ref{avel}) gives that $\sqrt{\frac{|\J_{n}|}{[\alpha n]}} \rightarrow 1$, as $n \rightarrow \infty$, a.s., therefore, by Lemma \ref{contog}, (\ref{sluts}), we have that 
\begin{equation}\label{eqaz}
\mathcal{L}\left(\frac{\sum_{x = - \alpha n}^{ \alpha n} \hat{\xi}_{2n}^{2\Z}(x)}{\sqrt{|\J_{n}|}}\right)  \xrightarrow{w} \mathcal{N}(0, \sigma^{2}),\mbox{ as } n \rightarrow \infty. 
\end{equation}
Hence, if we assume that 
\begin{equation}\label{propa}
\frac{\sum_{x \in \J_{n}} \hat{\xi}_{n}^{2\Z}(x)- \sum_{x = - \alpha n}^{ \alpha n} \hat{\xi}_{n}^{2\Z}(x)}{\sqrt{|\J_{n}|}}  \xrightarrow{p}0, \mbox{ as } n \rightarrow \infty, 
\end{equation}
then, (\ref{eqaz}) together with an application of Lemma \ref{contog}, (\ref{cteq1}), yields (\ref{eq12}), and the result is proved.  

We prove the remaining (\ref{propa}). To do this, we will show that the hypotheses of the general Proposition \ref{AnscC} are fulfilled when setting $(\hat{\xi}_{n}^{2\Z}(x),  l_{n}^{+}, r_{n}^{-})$ equal to $(X_{t}(j), m_{t}, M_{t})$ there. We have that:  
\begin{enumerate*}[label=\bf \itshape\alph*\upshape)] \item Recall that $\E(\xi_{n}^{2\Z}(x)) =  \rho_{n}$, where this equality comes from self-duality, see (\ref{eqduality}). We therefore have that assumption (\ref{mnzero}) holds since $(\hat{\xi}_{n}^{2\Z}(x))$ are centered r.v.'s.  
\item We now show that assumption (\ref{Cassass})  is granted for $(\hat{\xi}_{n}^{2\Z}(x))$ as follows. Note that, due to a corollary to Harris' correlation inequality [\cite{H77}], see [Thm.\ 2.14, Chpt.\ II; \cite{L85}], which applies since every deterministic configuration is positively correlated, we have that $(\xi_{n}^{2 \Z})$ has positive correlations for all $n$.  Because $\xi_{n}^{2 \Z}$ takes values on a partially ordered set, this gives that $\{\xi_{n}^{2 \Z}(x)\}$ are associated, and hence also $(\hat{\xi}_{n}^{2\Z}(x))$ are associated, because increasing functions of associated r.v.'s are also associated by using the definition.  \item   Because $\xi_{n}^{2 \Z}(x) \in \{0,1\}$, we have that $\E |\xi_{n}^{2 \Z}(x)| ^{2} \leq 1$, and therefore assumption (\ref{Cans0}) is also fulfilled. \item Furthermore, we have that (\ref{expdeccor}) gives that (\ref{Cansa}) is valid, because, by using that the covariance is a linear operation in the one argument if the other is fixed, we then have that $C^{-} = C^{+} = \frac{C}{1-\gamma} < \infty$. \item Finally, we have that (\ref{Cansb}) is valid for $m_{t} = l_{n}^{+}$ and $M_{t} = r^{-}_{n}$ due to (\ref{avel}). Hence, we have that (\ref{propa}) holds, and the proof is thus complete. 
\end{enumerate*}
\end{proof}

\begin{proof}[Proof of Proposition \ref{lemfindis}.] 
We will show that 
\begin{equation}\label{propeq1}
 \pr(A_{n} \geq a | \Omega_{\infty}) = \pr(\bar{A}_{n+k} \geq a | \tau = k), 
\end{equation}
for all $k \geq 0, a\in \R$. To see that proving the above display suffices complete this proof note that, due to Lemma \ref{lemtaf}, we have by the law of total probability and (\ref{propeq1}), that 
\begin{eqnarray*}
\pr(\bar{A}_{n + \tau} \geq a) &=& \sum_{k = 0}^{\infty} 
\pr(\bar{A}_{n + k} \geq a | \tau = k)\pr(\tau = k) \nonumber \\
& = & \pr(A_{n} \geq a | \Omega_{\infty}). 
\end{eqnarray*}

\noindent To state the next auxiliary lemma we require, we let $\mathcal{F}_{n}$ denote the $\sigma$-algebra associated to the part of the construction of the processes with bonds the end-vertices of which have time-coordinate no greater than $n$. 

\begin{lemma}\label{leminsid}
We let  
\begin{equation}\label{Omegaxn}
\xi^{(x,n)}_{m} = \{y: (x,n) \rightarrow (y,m+n)\},   \mbox{ } \Omega_{(x,n)} = \{ \xi^{(x,n)}_{m} \not= \emptyset, \mbox{ for all } m \geq 0\},  
\end{equation}
$ (x,n) \in \Lsite$,  $m \geq 0$. We have the following representation
\begin{equation}\label{repres}
\{ \tau = n, r^{-}_{n} = l^{+}_{n} = x \}  = \Omega_{(x, n)} \cap F,
\end{equation}
for some $F \in \mathcal{F}_{n-1}$.
\end{lemma}
\noindent We state next another general auxiliary statement we will use below. 
\begin{lemma}\label{couplemeta} 
Let $\eta^{-}, \eta^{+}$ be such that $\eta^{-}(O) =\eta^{+}(O) =1$ and $\eta^{-}(x)=0$, $ x\geq 2$, $\eta^{+}(x) = 0$, $x \leq -2$. Let $r_{n}^{\eta^{-}} = \sup\{x:  \xi_{n}^{\eta^{-}}(x) =1 \}$ and $l_{n}^{\eta^{+}}= \inf\{x: \xi_{n}^{\eta^{+}}(x) = 1\}$. We have that $\Omega_{\infty} =\{r_{n}^{\eta^{-}}  \geq l_{n}^{\eta^{+}}, \mbox{ for all } n \geq 1\}$. 
\end{lemma}

We now have that, for any $(x,k) \in \Lsite$,  
\begin{eqnarray}
 \pr(\bar{A}_{n+k} \geq a | \tau = k, r_{k}^{-} =  l^{+}_{k} = x)  &=&  \pr(\bar{A}_{n+k} \geq a | \Omega_{(x, k)}, r_{k}^{-} =  l^{+}_{k} = x, F) \label{aba1} \\ 
 & =&  \pr(\bar{A}_{n+k} \geq a | \Omega_{(x, k)}, r_{k}^{-} =  l^{+}_{k} = x)  \label{aba2}  \\
 & =&  \pr(\bar{A}_{n} \geq a |  \Omega_{\infty})  \label{aba3} \\ 
& =&  \pr(A_{n} \geq a | \Omega_{\infty}),  \label{aba4} 
\end{eqnarray}
where in (\ref{aba1}) we plug in  (\ref{repres}) from Lemma \ref{leminsid}, in (\ref{aba2}) we use independence of events measurable with respect to disjoint parts of $\mathscr{L}$ by construction, in (\ref{aba3}) we use translation-invariance with respect to $(x,k)$, and finally in (\ref{aba4}) we use that, by (\ref{basiccoup}),
$\bar{A}_{n} = A_{n}$ a.s.\ on $\Omega_{\infty}$.

The law of total probability gives
\begin{eqnarray}
\pr(\bar{A}_{n+k} \geq a | \tau = k)  & =& \sum_{x: (x, k) \in \Lsite} \pr(\bar{A}_{n+k} \geq a | \tau = k, r_{\tau}^{-} = x) \pr(r_{\tau}^{-} = x |\tau = k) \nonumber \\ 
& = & \pr(A_{n} \geq a | \Omega_{\infty})  \sum_{x : (x, k) \in \Lsite} \pr(r_{\tau}^{-} = x |\tau = k) \label{finiter} \\
& = & \pr(A_{n} \geq a | \Omega_{\infty}), \label{finiter2} 
\end{eqnarray}
where (\ref{finiter}) follows from (\ref{aba4}), and (\ref{finiter2}) follows from that, $\pr(|r_{k}^{-}| < \infty| \tau = k) = 1$, for all $k$, due to that $\pr(|r_{n}^{-}| < \infty) = 1$,  which follows from (\ref{avel}) by considering the contrapositive statement. This proof is thus complete.  
\end{proof}

We prove the remaining Lemmas \ref{lemtaf}, 
\ref{leminsid}, and \ref{couplemeta} that we stated and used above. 
 
\begin{proof}[Proof of Lemma \ref{lemtaf}.]
We will derive the estimate that there are $C, \gamma \in (0, \infty)$ such that
\begin{equation}\label{eqexpb}
\pr( \tau \geq n) \leq C e^{-\gamma n}, 
\end{equation}
for all $n \geq 1$.

Let $E^{r}_{n} = \{\exists m\geq n: r^{-}_{m} < 0\} \mbox{ and } E^{l}_{n} = \{\exists m\geq n: l^{+}_{m} >0\}.$ From (\ref{expoldev}), we have that 
\begin{equation}\label{ldev}
\pr(E^{r}_{n})\leq Ce^{- \gamma n},
\end{equation}
$ n \geq 1$. Further, note that 
\begin{equation}\label{coupltau}
\{ \tau \geq n\} \subseteq E_{n}^{r} \cup E_{n}^{l}, 
\end{equation}
where (\ref{coupltau}) follows by (\ref{taf}) and considering the contra-positive relation, i.e.\ that 
\[
\bar{E_{n}^{r}} \cap \bar{E_{n}^{l}} \subseteq  \{ \tau \leq n\}. 
\]
Hence, subadditivity and noting that by symmetry $\pr(E^{r}_{n}) = \pr(E^{l}_{n})$, gives
\[
\pr( \tau \geq n) \leq 2\pr(E^{r}_{n}),   
\]
from which the proof of (\ref{eqexpb}) is complete by (\ref{ldev}). 
\end{proof}

\begin{proof}[Proof of Lemma \ref{leminsid}.]
To prove (\ref{repres}),  note that it suffices to show that
\begin{equation}\label{rexp}
\tau =  \inf\{n \geq0: \cup_{(x,n) \in \Lsite}\{r^{-}_{n}= l^{+}_{n} = x\} \cap \Omega_{(x, n)}\}. 
\end{equation} 
However, Lemma \ref{couplemeta} and translation invariance give that, for any $(x, n) \in \Lsite$, 
\[
\Omega_{(x,n)} =  \{ r^{-}_{m} \geq l^{+}_{m}, \mbox{ } \forall m> n\}, \hspace{3mm} \mbox{ on } \{r^{-}_{n}= l^{+}_{n} = x\},
\]
hence (\ref{rexp})  is identical to (\ref{taf}), and (\ref{repres}) follows. 
\end{proof}

\begin{proof}[Proof of Lemma \ref{couplemeta}.]
Note that, since $\Omega_{\infty} = \cap_{n \geq 1} \Omega_{n}$, this statement follows directly from (\ref{rlc}).
\end{proof}

\begin{remark}\label{another}  \textup{ To derive the IP corresponding to Proposition \ref{bypro}, that is that, if we let $V_{t}^{n} = \frac{1}{\sigma \sqrt{\alpha n}}\sum_{x = - \alpha n}^{\alpha n} \hat{\xi}_{n}^{2\Z}(x)$, then we have that
\begin{equation}\label{eqWie}
V^{n} \Rightarrow W, 
\end{equation}
which means that the random functions $V^{n}$ converge to the Wiener measure, $W$, in the space $D$, see Section 13 in \cite{B68} for definitions and background in this regard. To derive (\ref{eqWie}) one may either modify the proof of Theorem 3 in \cite{NW81} to the case of arrays of r.v.'s in which the length of each row grows linearly with the row number, or alternatively, one may modify  the proof of Proposition \ref{bypro} in the
Appendix \ref{secApp} by invoking the IP extension of Lemma \ref{lem0} instead there. Letting $U_{t}^{n} = \frac{1}{\sigma \sqrt{\alpha n}}\sum_{x \in \J_{n}} \hat{\xi}_{n}^{2\Z}(x)$, we then have that by (\ref{eqWie}) and (\ref{avel}) the assumptions of 
Theorem 14.4 in \cite{B68} are fulfilled, hence yielding that 
\begin{equation}\label{Wie2}
U^{n} \Rightarrow W. 
\end{equation}
By appropriately modifying the argument from (\ref{taf}) onwards in the proof above we have that the corresponding IP to Theorem \ref{CLT} may also be derived from (\ref{Wie2}). Our proof of Theorem \ref{CLT} contrasts to the approach we briefly sketched here in that the former does not prerequisite
invoking any general statements. Further, note that the latter approach does not go through Proposition \ref{AnscC} and hence, does not provide with the random-indices CLT's we provide with in $\mathsection$ \ref{seccgprod}.} 
\end{remark}

\section{Proposition  \ref{AnscC} and Corollary \ref{cgprod}}\label{randev}

\subsection{Proof of Proposition  \ref{AnscC}.}\label{subspro}
The proof of Proposition \ref{AnscC} is divided into two parts. We will first derive Proposition \ref{AnscC} by means of relying on Lemma \ref{Ansc},  stated below here next, and proved below immediately thereafter, in this section. Prior to that, we 
also state here the Hajek-R\' enyi inequality for associated r.v.'s, due to \cite{C00}, see also \cite{S08}. Recall the definition of association given in $\mathsection \mathsection$ \ref{seccgprod}.  
\begin{lemma}\label{chr00}
Let $(X_{j}: j = 1, \dots, n)$ be associated r.v.'s such that $\E(X_{j})= 0$ for all $j$, and let also $(c_{j}: j = 1, \dots, n)$ be a sequence of non-increasing and positive numbers. Let $S_{k} = \sum_{j=1}^{k} X_{j}$. Then, we have that
\[
\pr\left(\max_{1 \leq k \leq n} c_{k} \left|S_{k}\right| \geq \epsilon \right) \leq 2\epsilon^{-2}\left(2 \sum_{j=1}^{n} c_{j}^{2} \Cov\left(X_{j}, S_{j-1} \right) + \sum_{j=1}^{n} c_{j}^{2} \E(X_{j}^{2})\right). 
\]
\end{lemma}

\noindent We may now state the following Lemma. 
\begin{lemma}\label{Ansc}
Let  $\{X_{t}(j):  j,  t \geq 0 \}$ be such that $\E(X_{t}(j))= 0$, $\forall$ $t, j$.  Let also ${S_{t}(i) = \sum_{j =0}^{i} X_{t}(j)}$. We assume the following: 
\begin{equation}\label{assass}
\{X_{t}(j):  j \geq 0 \} \mbox{ is associated for each } t. 
\end{equation}
\begin{equation}\label{ans0}
\sup_{j, t}\E(X_{t}(j)^{2})< \infty, 
\end{equation}
\begin{equation}\label{ansa}
\sup_{j, t}\Cov\left(X_{t}(j), S_{t}(j-1) \right) <\infty,
\end{equation}
Furthermore, we let $(N_{t}: t\geq 0)$ be integer-valued and non-negative r.v.'s, such that, for some $0<\theta<\infty$,
\begin{equation}\label{ansb}
\Law\left(\frac{N_{t}}{t}\right) \xrightarrow{w} \theta, \mbox{ as } t \rightarrow \infty. 
\end{equation}
Then, we have that 
\[
{\frac{S_{t}(N_t) - S_{t}([\theta t])}{ \sqrt{ [\theta t]} } \xrightarrow{p} 0,  \mbox{ as } t \rightarrow \infty.}
\]
\end{lemma}

\begin{proof}[Proof of Proposition \ref{AnscC}]
Let $M'_{t}= M_{t} \cdot 1\{M_{t} \geq 0\}$ and $m_{t}' = m_{t}\cdot 1\{m_{t} \leq 0\}$, where we recall that $M_{t}$ and $m_{t}$ are as in (\ref{Cansb}).
Note that in the notation introduced we have that 
\[
S_{t}(m'_{t}, M'_{t}) = S_{t}^{-}(m'_{t}) +S_{t}^{+}(M'_{t}), \mbox{ and } S_{t}(-[\theta t], [\theta t]) =  S_{t}^{-}([\theta t]) +S_{t}^{+}([\theta t]),  
\]
and thus,  by the triangle inequality, we have that 
\begin{equation}\label{trian}
\frac{|S_{t}(m'_t, M'_{t}) - S_{t}(-[\theta t], [\theta t])|}{ \sqrt{ [\theta t]}} \leq \frac{|S_{t}^{-}(m'_{t})-   S_{t}^{-}([\theta t])|}{ \sqrt{ [\theta t]}} + \frac{|S_{t}^{+}(M'_{t})-   S_{t}^{+}([\theta t])|}{ \sqrt{ [\theta t]}}. 
\end{equation}
However, the assumptions of Lemma \ref{Ansc} are appropriately satisfied, yielding that 
\[
\frac{|S_{t}^{-}(m'_{t})-   S_{t}^{-}([\theta t])|}{ \sqrt{ [\frac{\theta t}{2}]}} \xrightarrow{p} 0, \mbox{ and  }  \frac{|S_{t}^{+}(M'_{t})-   S_{t}^{+}([\theta t])|}{ \sqrt{ [\frac{\theta t}{2}]}} \xrightarrow{p} 0.
\]
as $t \rightarrow \infty$. Hence, from (\ref{trian}) and the display above,   
we have that 
\begin{equation}\label{SSzero2}
\frac{S_{t}(m_t', M_{t}') - S_{t}(-[\theta t], [\theta t])}{ \sqrt{ [\theta t]} } \xrightarrow{p} 0,  \mbox{ as } t \rightarrow \infty.
\end{equation}
To conclude the proof of (\ref{SSzero}), note that, again by the triangle inequality, 
\begin{eqnarray*}\label{trianagain} 
\frac{|S_{t}(m_t, M_{t}) - S_{t}(-[\theta t], [\theta t])|}{ \sqrt{ [\theta t]} } &\leq & \frac{|S_{t}(m_t, M_{t}) - S_{t}(m_t', M_{t}')|}{ \sqrt{ [\theta t]} } +  \nonumber \\
&+& \frac{|S_{t}(m_t', M_{t}') - S_{t}(-[\theta t], [\theta t])|}{ \sqrt{ [\theta t]}},
\end{eqnarray*}
so that in view of (\ref{SSzero2}), it suffices to show that 
\begin{equation}\label{deslem}
\lim_{t \rightarrow \infty} \pr(|S_{t}(m_t, M_{t}) - S_{t}(m'_t, M'_{t})| > \epsilon) = 0,
\end{equation}
which follows simply by noting that, for all $\epsilon >0$,  
\[
\pr(|S_{t}(m_t, M_{t}) - S_{t}(m'_t, M'_{t})| > \epsilon) \leq \pr(M_{t} \leq -1 \mbox{ or } m_{t} \geq 1),  
\]
and hence \ref{deslem} follows, since from (\ref{Cansb}) we have that $\lim_{t \rightarrow \infty} \pr(M_{t} \leq -1) = 0$ and $\lim_{t \rightarrow \infty} \pr(m_{t} \geq 1) = 0$. 
The proof is thus complete. 
\end{proof}

\begin{proof}[Proof of Lemma \ref{Ansc}]

Note that without loss of generality we may take $\theta = 1$. Let $\epsilon \in (0,1)$, and also let $m(t) = [t(1-\epsilon^{3})]+1$ and $n(t) =  [t(1+\epsilon^{3})]$. Let $Y_{i}(t) = X_{t}(t+i)$, $i = 1, \dots, n(t)$ and denote their partial sums as $Z_{k}(t) = \sum_{i=1}^{k} Y_{i}(t)$,  then we have that 
\begin{eqnarray}
\max_{k=t, \dots, n(t)} |S_{t}(k) - S_{t}(t)|   &= &\max_{k=t+1, \dots, n(t)} \left|\sum_{j = t+1}^{k} X_{t}(j) \right| \nonumber \\
&=& \max_{k=1, \dots, [t\epsilon^{3}]} |Z_{k}(t)|. 
\end{eqnarray}
From (\ref{assass}) we have that Lemma \ref{chr00} applies and choosing there $c_{k} = \frac{1}{\sqrt{t}}$,  gives that 
\begin{eqnarray}\label{above0}
\pr \left(\max_{k=t, \dots, n(t)} |S_{t}(k) - S_{t}(t)|   \geq \epsilon \sqrt{t}\right) &=& \pr\left( \max_{k=1, \dots, [t\epsilon^{3}]} |Z_{k}(t)|   \geq \epsilon \sqrt{t}\right) \nonumber \\
&\leq & \frac{2}{\epsilon^{2} t}  \left(2 \sum_{j =1}^{[t \epsilon^{3}]} \Cov(Y_{j}(t), Z_{j-1}(t)) + \sum_{j =1}^{[t \epsilon^{3}]} \E(Y_{j}(t))^{2} \right)  \nonumber \\
&\leq &  C \epsilon, \label{cs} 
\end{eqnarray}
where  $C$ is independent of $t$, and  (\ref{cs}) comes from (\ref{ans0}) and (\ref{ansa}). Similarly, letting $Y'_{i}(t) =X_{t}(t+1-i)$, for $i = 1, \dots t - m(t) +1$ and $Z'_{k}(t) = \sum_{j=1}^{k}Y'_{i}(t)$,  we have that
\begin{eqnarray}
\max_{k=m(t), \dots, t} |S_{t}(k) - S_{t}(t)|   &= &\max_{k=m(t), \dots, t-1} \left|\sum_{j = k+1}^{t} X_{t}(j) \right|  \nonumber \\
& = & \max_{k = 1, \dots, [t \epsilon^{3}]-1} |Z_{k}'(t)|. 
\end{eqnarray}
Again, we can apply Lemma \ref{chr00} with $c_{k} = \frac{1}{\sqrt{t}}$ from (\ref{assass}), so that 
\begin{eqnarray}\label{above}
\pr \left(\max_{k=m(t), \dots, t} |S_{t}(k) - S_{t}(t)| \geq \epsilon \sqrt{t}\right) &=& \pr\left(\max_{k = 1, \dots, [t \epsilon^{3}]-1} |Z_{k}'(t)| \geq \epsilon \sqrt{t}\right) \nonumber \\
&\leq & \frac{2}{\epsilon^{2} t}  \left(2 \sum_{j =1}^{[t \epsilon^{3}]-1} \Cov(Y'_{j}(t), Z'_{j-1}(t)) + \sum_{j =1}^{[t \epsilon^{3}]-1} \E(Y_{j}'(t))^{2} \right)  \nonumber \\
&\leq &  C \epsilon, \label{ulti2}
\end{eqnarray}
where again we use (\ref{ans0}) and (\ref{ansa}) in (\ref{ulti2}). 

Partitioning according to the event $N_{t} \in [m(t), n(t)]$ and its complement and then using (\ref{above0}) and (\ref{above}) gives that
\begin{eqnarray}\label{lastpro}
\pr(|S_{N_{t}} - S_{t}| \geq \epsilon \sqrt{ t}) &\leq& \pr(|S_{N_{t}} - S_{t}| \geq \epsilon \sqrt{ t}, N_{t} \in [m(t), n(t)]) + \pr(N_{t} \not\in [m(t), n(t)])  \nonumber \\ 
& \leq &  \pr \left( \max_{m(t)\leq k \leq t} |S_{k} - S_{t}| \geq \epsilon \sqrt{t}\right)+ \nonumber \\ 
&& + \hspace{2mm} \pr\left(\max_{t \leq k \leq n(t)} |S_{k} - S_{t}| \geq \epsilon \sqrt{ t}\right)+ \pr(N_{t} \not\in [m(t), n(t)]) \nonumber  \\
&\leq& 2C \epsilon  + \pr\left(N_{t} \not\in [m(t), n(t)]\right),
\end{eqnarray}
where for the last inequality we invoke (\ref{cs}) and (\ref{above}). However,  (\ref{ansb}) gives
\[
\limsup_{t\rightarrow \infty}\pr\left(N_{t} \not\in [m(t), n(t)]\right) = 0, 
\]
and hence, from (\ref{lastpro}) we get that 
\[
\limsup_{t\rightarrow \infty}\pr(S_{N_{t}} - S_{t} \geq \epsilon \sqrt{ t}) \leq 2C \epsilon,
\]
which due to that $\epsilon$ is arbitrary, completes the proof. 
\end{proof} 
\vspace{2mm}

\subsection{Proof of Corollary \ref{cgprod}}\label{prfcor}

The following Theorem, due to \cite{CG84}, is applied in the proof of this Corollary following next. 
\begin{theorem}\label{TCG}
Let $\{X_{n}(j): 0 \leq j \leq n\}$ be such that $\E(X_{n}(j))= 0$, $\forall$ $n, j$, and suppose that  (\ref{eqcltassoc}), (\ref{cg1}), and  
(\ref{cg2}) hold.  Letting $S_{n}(n) = \sum_{j =0}^{n} X_{n}(j)$, we then have that
\begin{equation}\label{eqcgt}
\frac{S_{n}(n)}{\sqrt{n}} \xrightarrow{w} \mathcal{N}(0,\sigma^{2}),  \mbox{ as } n \rightarrow \infty, 
\end{equation}
where $\sigma^{2}:= \lim_{n \rightarrow \infty} \textup{\mbox{Var}} (S_{n}(n) / \sqrt{n})$, $0<\sigma^{2}<\infty$. 
\end{theorem}

\begin{proof}[Proof of Corollary \ref{cgprod}.]
Note that it suffices to only show the first conclusion Corollary \ref{cgprod}, for the second one follows from that by an application of Lemma \ref{contog}. Note also that there is no loss of generality in assuming $\theta =1$. 
We thus have that the hypotheses of Theorem \ref{TCG} are met, and hence, (\ref{eqcgt}) holds. From it and Lemma \ref{Ansc}, the proof is complete by an application of Lemma \ref{contog}.   
\end{proof}

\section{Appendix}\label{secApp}

\begin{proof}[Proof of Proposition \ref{bypro}.] 
We let $(\xi_{n}^{2\Z})$ and $(\bar{\xi}_{n})$ be the processes with starting sets $2\Z$ and $\bar{\xi}_{0}\sim \bar{\nu}$, where we enlarge our probability space to support a $\bar{\nu}$-distributed independent random set ${\bf S}$ by setting $\bar{\xi}_{0} = S$, on $\{{\bf S} = S\}$. We also let $K_{n}(a)$ be the set of points of $\Lsite$ inside a cone of slope $a>0$ and apex $(O,0)$, as follows $K_{n}(a) = \{x: -an \leq x \leq an \mbox{ and } (x,n) \in \Lsite\}$, $n \geq 0$. We have the next Lemma. 

\begin{lemma}\label{lema1} $\{\forall x \in K_{n}(a), \xi_{n}^{2 \Z}(x) =  \bar{\xi}_{n}(x)\}, \mbox{ } \forall n \mbox{ large, $\pr$-a.s.}$
\end{lemma}

\begin{proof}
From the Borel-Cantelli lemma and the union bound, since $|K_{n}(a)|$ grows linearly in $n$,  it suffices to show that, for all $a>0$, there are $C, \gamma$  such that, for any $x \in K_{n}(a)$,
\begin{equation*}\label{desirag}
\pr(\bar{\xi}_{n}(x) \not= \xi^{2 \Z}_{n}(x)) \leq  C e^{- \gamma n}, 
\end{equation*}
$n \geq 1$. However, we have that 
\begin{eqnarray*}
\pr(\bar{\xi}_{n}(x) \not= \xi^{2 \Z}_{n}(x)) &=&  \pr(\xi^{2 \Z}_{n}(x) = 1, \bar{\xi}_{n}(x) = 0) \\
&=&  \pr(\xi^{2 \Z}_{n}(x) =1 ) - \pr(\bar{\xi}_{n}(x) =1 )  \label{eqone} \\
 &=& \pr(\Omega_{n}) - \pr(\Omega_{\infty})  \label{eqone2*} \\ 
& \leq &  C e^{- \gamma n}, \label{eqone4}  
\end{eqnarray*}
$n \geq 1$, where in the first line we used that, by (\ref{eqmon}), $\bar{\xi}_{n} \supseteq \xi^{2 \Z}_{n}$, and in the second one we used that and, in addition, the law of total probability; in the third line we used (\ref{eqduality}), and also (\ref{duapp2}) together with stationarity; and the last inequality comes from (\ref{lemexp}). This completes the proof. 
\end{proof} \vspace{-3mm}
\noindent We also note here the following consequence of Lemma \ref{lem0} above. 
\begin{corollary}\label{cor0} $\displaystyle{\mathcal{L}\left(
\frac{\sum_{x \in K_{n}(\alpha)} (\bar{\xi}_{n}(x) - \rho)}{\sqrt{\alpha n}}\right) \xrightarrow{w} \mathcal{N}(0, \sigma^{2})}$, as $n \rightarrow \infty$. 
\end{corollary}

\begin{proof}
Let $\xi_{n}^{'}$ be such that 
\begin{equation*}\label{xistar}
  \xi_{n}^{'}(x) = 
    \begin{cases}
      \bar{\xi}_{n}(x), & \text{if} \ n \in 2\Z_{+}  \\
       \bar{\xi}_{n}(x-1) &  \text{if} \ n \in 2\Z_{+}+1, 
    \end{cases}
  \end{equation*}
and note that, for all $n$, $\xi_{n}^{'} \sim \bar{\nu}$. Applying Lemma \ref{lem0}, completes the proof. 
\end{proof}
\vspace{-3mm}
\noindent Note that Lemma \ref{lema1} gives that
$\sum_{x \in K_{n}(\alpha)} \xi^{2 \Z}_{n}(x)   = \sum_{x \in K_{n}(\alpha)}\bar{\xi}_{n}(x),   \forall n \mbox{ large,  $\pr$-a.s.}$, so that from Corollary \ref{cor0}, and Lemma \ref{contog}, (\ref{cteq1}), we have that
\begin{equation*}\label{dedeq}
\mathcal{L}\left(
\frac{\sum_{x \in K_{n}(\alpha)} (\xi^{2 \Z}_{n}(x) - \rho)}{\sqrt{\alpha n}}\right) \xrightarrow{w} \mathcal{N}(0, \sigma^{2}), \mbox{ as } n \rightarrow \infty. 
\end{equation*}
By the last display above, since $\hat{\xi}_{n}^{2\Z}(x) = \xi_{n}^{2\Z}(x) - \rho_{n}$, again by applying Lemma \ref{contog}, (\ref{cteq1}), we have that it suffices to show that $ \displaystyle{\sqrt{\alpha n} (\rho_{n} - \rho) \rightarrow 0}$, as $n \rightarrow \infty$, which holds by (\ref{lemexp}), and hence the proof is complete.
\end{proof}

\noindent\textbf{Acknowledgments:} 
This work has been supported during non-overlapping periods of time by CONICET, by FAPESP grant 2016/03988-5, and, currently, by PNPD/CAPES.

\bibliographystyle{plain}

\vspace{-10mm}
\textsc{ \\
\noindent \footnotesize Postal address: \\ Instituto de Matem\' atica e Estat\' istica,\\
Universidade de S\~ ao Paulo\\
Rua do Mat\~ ao, 1010\\
CEP 05508-900- S\~ ao Paulo\\
Brasil}

\end{document}